\documentclass{ws-procs9x6}

\numberwithin{equation}{section}
\def\phi{\varphi }
\def\reel{{\mathbb R}}
\def\nat{{\mathbb N}}

\def\comp{{\mathbb C}}
\def\ganz{{\mathbb Z}}

\def\bar{\overline }
\newtheorem{randomwalks}[theorem]{{\bf Random walks}}

\newtheorem{ex-dunklhypergroup}[theorem]{{\bf Examples: Hypergroups on Weyl chambers
  associated with Dunkl operators}}
\newtheorem{ex-complexliegroup}[theorem]{{\bf Examples: Hypergroups associated with
complex semisimple Lie groups with finite center}}

\begin{document}
\title{Deformations of convolution semigroups on commutative hypergroups  }
\author{MARGIT R\"OSLER}
\address{Mathematisches Institut, Universit\"at G\"ottingen\\ 
Bunsenstr. 3--5, D-37073 G\"ottingen,
Germany\\
E-mail: roesler@uni-math.gwdg.de}
\author{ Michael Voit}
\address{ Fachbereich Mathematik,  Universit\"at Dortmund\\
D-44221 Dortmund, Germany\\
E-mail: michael.voit@mathematik.uni-dortmund.de}

\maketitle

%\noindent
% {\bf Mathematics Subject  Classification:}
%{Primary:~43A10;
%Secondary:~43A62, 60B15, 47D07.}

%\noindent
% {\bf Keywords:}
%{Gelfand pairs, deformation of hypergroups,
%positive semicharacters, convolution semigroups, L\'evy measures.}

\abstracts{
It was recently shown by the authors that deformations
 of hypergroup convolutions w.r.t.
positive semicharacters can be used to explain probabilistic 
connections between 
the Gelfand pairs $(SL(d,\comp), SU(d))$ and Hermitian matrices. 
We here study connections between general convolution semigroups on 
commutative hypergroups and their deformations.
 We are able to develop a satisfying theory, if
the underlying positive semicharacter has some growth property. 
We present several examples which indicate that this  growth condition
holds in many  interesting cases. 
}

\section{Introduction}
 Klyachko \cite{Kl2} recently derived a connection between
$SU(d)$-biinvariant  random walks  on  $SL(d,\comp)$
and  
random walks on the  additive group ${\bf H}_{d,0}$ of 
all hermitian $d\times d$-matrices with trace 0,
 whose transition probabilities are
invariant under 
conjugation by  $SU(d)$. He  used 
this connection to transfer the recent solution of the spectral problem for sums of 
hermitian matrices (\cite{Kl1}, \cite{KT}) 
to  the possible singular spectrum of products of
 random matrices from  $SL(d,\comp)$ with given singular spectra. The
 singular spectrum of a matrix $A\in SL(d,\comp)$ here means the spectrum of
 the positive definite matrix $\sqrt{AA^*}$. 
Klyachko's  connection between  $SL(d,\comp)$ and  ${\bf H}_{d,0}$ was  
explained in a different way and extended  by the authors in \cite{RV};
it is shown in \cite{RV} that the commutative Banach algebra
of all $SU(d)$-biinvariant bounded 
 measures on  $SL(d,\comp)$ may be embedded into 
the Banach algebra of all bounded 
measures on the Euclidean space  ${\bf H}_{d,0}$
in an isometric, probability preserving way. The proof of this fact, 
which has some applications in probability theory (see \cite{RV}), 
depends on so-called
  deformations of  hypergroup convolutions with respect to
 positive semicharacters as introduced in \cite{V1}.
These    deformations  lead
 to  
connections between random walks
and convolution semigroups on different, but closely related hypergroups.
This forms 
the motivation to  investigate systematically
when and how convolution semigroups of probability measures
on a commutative hypergroup $(X,*)$ can be transformed canonically
into convolution semigroups on a   deformation $(X,\bullet)$ of $(X,*)$.
In particular we show that the generators and L\'evy measures 
of the original and the deformed  convolution semigroup are closely related
 whenever this  transformation is possible. We mention that the
deformation of  convolution semigroups is closely related to
 Doob's $h$-transform, and that  L\'evy processes
associated with a convolution semigroup and its deformation are 
related by a Girsanov transformation on the path space; see \cite{V2}.

The paper is organized as follows: In Section 2 we collect  some 
 facts on deformations  and present  examples.
In particular we indicate how for a maximal compact subgroup $H$ of
a complex, non-compact, connected semisimple Lie group  $G$,
the double coset hypergroup $G//H$ may be regarded as  deformation 
of an orbit hypergroup. This includes 
the examples above.  Section 3 is  devoted to
deformations of  convolution semigroups w.r.t.  positive semicharacters
$\alpha_0$. We show that this concept works in a satisfying way 
 under a canonical growth
 condition on the  convolution semigroup together with some growth
condition concerning $\alpha_0$.  Section 4 finally contains
examples  where this 
 condition on $\alpha_0$ is satisfied.
In fact, we have no  example for which this condition would not hold.

 \section{Deformations of commutative hypergroups}

We  give a quick introduction. First,
let us  fix   notations.
For a locally compact Hausdorff space  $X$, $M^+(X)$ denotes the 
 space of all positive 
Radon measures on $X$, and  $M_b(X)$ the Banach space of all bounded
regular complex Borel measures  with the total variation norm.
Moreover,  $M^1(X)\subset M_b(X) $ is the set of all probability measures,
$M_c(X)\subset M_b(X) $  the set of all measures with compact support, and 
 $\delta_x$   the point
measure in $x\in X$.
The spaces $C(X)\supset C_b(X)\supset C_0(X) \supset C_c(X)$ 
of continuous functions are given  as usual. 

\begin{definition} A hypergroup $(X,*)$ 
consists of a locally compact Hausdorff space $X$ and a weakly continuous,
  probability preserving convolution $*$
 on $M_b(X)$ 
such that $(M_b(X),*)$ is a Banach algebra and $*$ preserves compact
supports.
 Moreover, there  exists an identity $e\in X$ (such that  $\delta_e$
is the identity of   $(M_b(X),*)$) as well as  
  a continuous involution $x\mapsto\bar x$ on $X$ that replaces the group 
inverse. For  details  we refer to \cite{BH} and  \cite{J}.

We here only deal with commutative hypergroups  $(X,*)$, i.e., 
 $*$ is commutative. In this case 
there exists an (up to normalization)  unique  
Haar measure $\omega\in M^+(X)$ which is characterized by
$\, \omega(f)=\omega(f_x)$ for all  $\,f\in C_c(X)$ and 
$x\in X$, 
where we use the notation
\[f_x(y):=f(x*y):= \int_X f\> d(\delta_x*\delta_y).\]
Similar to the dual  of a locally compact abelian group, 
one defines
\begin{align}
\chi(X):=&\{\alpha\in C(X):\> \alpha\ne 0,\>\>  \alpha(x*y)=
 \alpha(x)\alpha(y)\>\>\text{for all}\>\> x,y\in X\},\notag\\
 X^*:=&\{\alpha\in\chi(X):\> \alpha(\bar x)=
\overline{\alpha(x)}\>\>\text{for}\>\> x\in X\};  \quad
\widehat X:=X^*\cap C_b(X).\notag
\end{align}
Elements of $X^*$ and  $\widehat X$ are called
semicharacters and characters respectively.
All spaces  are locally compact  w.r.t. 
 the compact-uniform  topology.  
\end{definition}

\begin{example}\label{ex_1}
\begin{enumerate}
\item[\rm{(1)}]   Let $K$ be a compact subgroup of a locally compact group $G$.
Then 
$$M_b(G\|K):=\{\mu\in M_b(G):\> \delta_x*\mu*\delta_y=\mu\>\>\text{for all}\>\> x,y\in K\}
$$
is a Banach-$*$-subalgebra of $M_b(G)$ with 
 the normalized Haar measure $dk\in M^1(G)$ of $K$ 
as identity.
The double coset space $G//K:=\{KxK:\> x\in G\}$ is   locally compact
 w.r.t. the quotient topology, and  the canonical projection 
$p:G\to G//K$ induces a  probability preserving, isometric isomorphism
$\,p: M_b(G\|K)\to M_b(G//K)$ of Banach spaces by taking 
images of measures. The transport of
the convolution on  $M_b(G\|K)$ to $M_b(G//K)$ via $p$ leads to a 
 hypergroup
structure $(G//K, *)$ with identity $K$ and involution
$(KxK)^-:=Kx^{-1}K$, and  
$p$  even becomes  a Banach-$*$-algebra isomorphism.
If $G//K$ is commutative, i.e., $(G,K)$
  is a Gelfand pair, then a $K$-biinvariant function $\phi\in C(G)$ 
with $\phi(e)=1$ is   spherical if 
$\phi(x)\phi(y)=\int_K f(xky)\> dk$  for $x,y\in G$.
The  functions $\alpha\in \chi(G//K)$ are in one-to-one
correspondence with the spherical functions on $G$ via 
 $\alpha\mapsto \alpha\circ p$ for the canonical projection $p:G\to G//K$.
\item[\rm{(2)}] Let $(V, \langle\,.\,,.\,\rangle)$ be a 
finite-dimensional Euclidean vector space and
 $K\subset O(V)$   a compact subgroup  of the orthogonal
group of $V$. For $\mu\in M_b(V)$, 
denote the image measure of $\mu$ under $k\in K$ by $k(\mu).$
Then  the space of $K$-invariant measures
$$M_b^K(V):=\{\mu\in M_b(V):\> k(\mu)=\mu\>\>\text{for all}\>\>
k\in K\}$$
is a Banach-$*$-subalgebra of 
 $M_b(V)$ (with the group convolution) with identity $\delta_0$.
The space $V^K:=\{K.x:\> x\in V\}$
 of all $K$-orbits in $V$ is  again locally compact,
 and   the canonical projection 
$p:V\to V^K$  induces a probability preserving,
 isometric isomorphism
$p: M_b^K(V)\to M_b(V^K)$ 
of Banach spaces and an  associated
 orbit hypergroup structure  $(V^K, *)$ such that
$p$  becomes  an isomorphism of Banach-$*$-algebras. The involution 
on $V^K$  is given by $\overline{K.x} = -K.x$. Moreover, the continuous
 functions
\begin{equation}\label{char-eucl}
 \alpha_\lambda(K.x) =
 \int_K e^{i\langle \lambda,\,k.x\rangle} dk \quad (x\in V)
\end{equation}
are   multiplicative on  
$(V^K,*)$ for  $\lambda\in V_\comp$, the complexification of $V$, and
$\alpha_\lambda \equiv \alpha_\mu$   if and only if
 $K.\lambda = K.\mu$. It is  known (see \cite{J}) that 
$\widehat{V^K}=\{\alpha_\lambda:\> \lambda\in V\}$. 
\end{enumerate}
\end{example}

 By  \cite{V1}, positive semicharacters lead to 
deformed  convolutions:

\begin{proposition} 
Let $\alpha_0\in X^*$
be a positive semicharacter on the commutative hypergroup $(X,*)$,
 i.e.,   $\alpha_0(x) >0$ for  $x\in X$.
Then 
\begin{equation}\label{deformconvo}
\mu\bullet\nu = \alpha_0((\alpha_0^{-1}\mu)*(\alpha_0^{-1}\nu))
\quad\quad\quad(\mu,\nu\in M_c(X))
\end{equation}
 extends uniquely  to a bilinear, associative,
probability preserving, 
weakly  continuous convolution $\bullet$ on $M_b(X)$, and
 $(X,\bullet)$
becomes 
a commutative hypergroup with the identity and involution of $(X,*)$.
 $(X,\bullet)$ will be called deformation of $(X,*)$
w.r.t. $\alpha_0$. 
\end{proposition}

 Eq.(\ref{deformconvo}) shows that  $\mu\mapsto \alpha_0\mu$ is 
 an  algebra isomorphism between $(M_c(X),*)$ and 
 $(M_c(X),\bullet)$ which for unbounded $\alpha_0$
cannot be extended to $M_b(X)$; cf. Section 3.

Many data of $(X,\bullet)$ can be expressed in
terms of $\alpha_0$ and corresponding  data of   $(X,*)$. For instance, if
 $\omega$ is a Haar measure of  $(X,*)$, then $\alpha_0^2\omega$
is a Haar measure of $(X,\bullet)$. Moreover,
the mapping $M_{\alpha_0}: 
  \,\alpha\mapsto \alpha/\alpha_0\,$ is a homeomorphism 
 between $(X,*)^*$ and $ (X,\bullet)^*$, and also
  between $\chi(X,*)$ and $\chi(X,\bullet)$; see \cite{V1} and  \cite{RV}.

\begin{remark} 
Deformation is  transitive as follows:
Let $(K,\bullet)$ be the deformation of  $(K,*)$ w.r.t. 
 $\alpha_0$, and let $\beta_0$ be a 
positive   semicharacter on  $(K,\bullet)$. Consider further the
deformation  $(K, \diamond)$ of $(K,\bullet)$ w.r.t. $\beta_0$. 
The function $\alpha_0\beta_0$ is a positive  semicharacter on  $(K,*)$,
and  $(K, \diamond)$ is the deformation  of  $(K,*)$ w.r.t.
  $\alpha_0\beta_0$. For $\beta_0=1/\alpha_0$, one
  obtains $\diamond=*$.
\end{remark}

We next present some  examples; for 
further examples  see
 Section 4. 

\begin{example}\label{deform-euc}
 Let $(V, \langle\,.\,,\,.\,\rangle)$ be an  $n$-dimensional  
Euclidean vector space,
 $K$ a  compact subgroup  of the orthogonal group $O(V)$, and
  $(V^K, *)$   the associated orbit hypergroup.
  Fix  $\rho\in V$ with $-\rho\in K.\rho$,
and consider the  function 
$e_\rho(x):=e^{\langle \rho,x\rangle}$
on $V$ and
$$M_c^{\rho,K}(V):=\{e_\rho\mu:\>
 \mu\in M_c(V) \>\>\text{$K$-invariant}\}.$$
The multiplicativity of  $e_\rho$  on $V$
yields that w.r.t. the  group convolution on  $M_c(V)$, we have
$ e_\rho\mu* e_\rho\nu=  e_\rho(\mu*\nu)$. Hence, 
$M_c^{\rho,K}(V)$ is a subalgebra of 
 $M_b(V)$, and  its norm-closure 
$$M_b^{\rho, K}(V):= \overline{M_c^{\rho,K}(V)}$$
 a Banach subalgebra.
On the other hand, $\displaystyle \alpha_0(K.x):=\int_K e_\rho(k.x)dk \, 
(x\in V)$  is a positive semicharacter on  $(V^K, *)$;
see  Example \ref{ex_1}(2) above as well as Proposition 2.8 of  \cite{RV}.
 Proposition 2.8 of  \cite{RV} also  states that for the deformation
 $(V^K,\bullet)$  of 
  $(V^K, *)$  w.r.t.  $\alpha_0$,  the canonical projection 
$p: V\to V^K$ induces  a
probability preserving isometric isomorphism of Banach algebras  
from  $M_b^{\rho, K}(V)$ onto  $M_b(V^K, \bullet)$.
In other words, the deformed hypergroup algebra may be regarded as 
 Banach algebra of (not longer $K$-biinvariant) measures on $V$.
\end{example}

\begin{example}
It is well-known that the double coset hypergroup $SL(2,\comp)//SU(2)$
and the   orbit hypergroup  $(\reel^3)^{SO(3)}$ may be 
identified with  $[0,\infty[$, and that the associated
hypergroup structures on   $[0,\infty[$ are deformations of each other; 
 see  \cite{BH}, \cite{RV}, or \cite{V1}.
\end{example}

Here is the higher rank extension of this example:

\begin{example}
Let $G$ be a complex, noncompact, connected
semisimple  Lie group with finite center
and $K$  a maximal compact subgroup.
Consider the   Cartan
decomposition $\mathfrak g= \mathfrak k+ \mathfrak p$ of the Lie algebra of
$G$, and choose a  maximal abelian subalgebra $\mathfrak a\subset \mathfrak p$.
 $K$ acts on $\mathfrak p$ via the adjoint representation as a
group of orthogonal transformations w.r.t. the Killing form
$\langle\,.\,,\,.\,\rangle$ as
scalar product. Let  $W$ be the Weyl group of $K$, which acts on
  $\mathfrak a$ as  finite reflection group;  here and further on we identify $\mathfrak a$
 with its dual $\mathfrak a^*$ via the Killing form. Fix some Weyl
  chamber $\mathfrak a_+\subset \mathfrak a$ and the associated set
 $\Sigma^+$ of positive roots.
 Then the closed chamber
$C:=\overline{\mathfrak a_+}$ is a fundamental domain for the action
of $W$ on $\mathfrak a$, and   $C$ can be   identified with 
 the orbit hypergroup $(\mathfrak p^K, *)$, where  a
$K$-orbit in $ \mathfrak p$ corresponds to its representative
in $C$. 

 $C$ can also  be identified
with the commutative double coset hypergroup
  $G//K$ where  $x\in C$ corresponds
to the double coset $K(e^x)K$. Denote the corresponding  
  convolution by $\bullet$.
Using  the known 
formulas for the spherical functions on $G//K$ and 
 $\mathfrak p^K$ (see Helgason \cite{He}), we proved
in  \cite{RV} that  $(G//K, \bullet)=(C,\bullet)$
is the deformation of the orbit hypergroup
 $(\mathfrak p^K, *)=(C,*)$ w.r.t.  the positive semicharacter
$\alpha_{-i\rho}$ (in the sense of Example \ref{ex_1}(2)) with
$$\textstyle
\rho:= \sum_{\alpha\in\Sigma^+}\alpha\in \mathfrak a_+\,.$$
As 
$-\rho\in K.\rho$,  the construction in Example 
\ref{deform-euc} shows that  $M_b(G||K)$ may be embedded
into $M_b(\mathfrak p)$ in an isometric, probability preserving way. 
Here are the most prominent examples (c.f.
Appendix C of \cite{Kn}).
\begin{enumerate}\itemsep=-1pt
\item[\rm{(1)}] {\bf The $A_{d-1}$-case.}\,
$K=SU(d)$ is a maximal compact subgroup of  $G=SL(d,\comp)$. 
In the   Cartan decomposition
$\mathfrak g= \mathfrak k+ \mathfrak p$ we obtain   $\mathfrak p$ as
 the additive group ${\bf H}_d^0$ of
all Hermitian $d\times d$-matrices with trace $0$, on which
$SU(d)$ acts by conjugation. Moreover, $\mathfrak a$
consists of all real diagonal matrices with trace 0 and will be identified
with
\[\{x = (x_1,\ldots,x_d)\in\reel^d:\> \sum_i x_i=0\}\]
on which the Weyl group acts as the symmetric group $S_d$.
We thus   take
 $$C:=\{x = (x_1,\ldots,x_d)\in\reel^d:
 \> x_1\ge x_2\ge\ldots\ge x_d,\> \sum_i x_i=0\}.$$ 
Then in particular,
$\,\rho= \left({d-1},{d-3}, \ldots,{-d+3},{-d+1}\right)$.
\item[\rm{(2)}] {\bf The $B_d$-case.} For $d\ge 2$ consider 
$G=SO(2d+1,\comp)$ with  maximal compact subgroup $K=SO(2d+1)$.
Here $\mathfrak a$ may be identified with $\reel^d$, and we may choose
$$C=\{x\in\reel^d:\> x_1\ge x_2\ge\cdots\ge x_d\ge0\}$$
with the Weyl group  $W\simeq S_d \ltimes\ganz_2^d$, and
$\rho=\left(2d-1,2d-3,\ldots,1\right).$
\item[\rm{(3)}] {\bf The $C_d$-case.} For $d\ge 3$ let
$G=Sp(d,\comp)$ with  the maximal compact subgroup $K=Sp(2d+1)$.
Here,  $\mathfrak a=\reel^d$ with 
 $C$ and  $W$ as in the $B_d$-case.
We  have
 $\rho=(2d,2d-1,\ldots,2)$.
The preceding results
on hypergroup deformations 
imply that the  hypergroups $Sp(d,\comp)//Sp(2d+1)$ and
$SO(2d+1,\comp)//SO(2d+1)$ are (up to isomorphism) 
deformations of each other; see also \cite{RV} .
\item[\rm{(4)}] {\bf The $D_d$-case.} For $d\ge 4$ let
$G=SO(2d,\comp)$ with maximal compact subgroup $K=SO(2d)$.
In this case  $\mathfrak a=\reel^d$ and we may take
$$C=\{x\in\reel^d:\> x_1\ge x_2\ge\cdots\ge x_{d-1}\ge |x_d|\}$$
with 
$\rho=(2d-2,2d-4,\ldots,2,0).$
\end{enumerate}
\end{example}

\section{Deformation  of  convolution semigroups}

We  now always assume  that $\alpha_0$ is a positive semicharacter 
on a  $\sigma$-compact, second countable
 commutative hypergroup $(X,*)$ and that  $(X,\bullet)$ 
is the associated deformation. We  show how under a natural growth
condition, 
convolution semigroups on $(X,*)$ can be deformed into 
convolution semigroups on  $(X,\bullet)$. To describe this 
 condition, we introduce the spaces
\begin{align}
&M_{\alpha_0}^{b,+}(X):= \,\{ \mu\in M^{b,+}(X):\> 
 \alpha_0\mu\in M^{b,+}(X)\},\notag\\
&M_{\alpha_0}^1(X):=\, M^1(X)\cap M_{\alpha_0}^{b,+}(X)\notag
\end{align}
as well as the transformation
$$R_{\alpha_0}: M_{\alpha_0}^{b,+}(X)\to  M^1(X), \quad
 \mu\mapsto \frac{1}{\int_X\alpha_0\> d\mu} \cdot \alpha_0\mu.$$
%which  satisfies $R_{\alpha_0}( M_{\alpha_0}^1(X))\subset  M^1(X)$.

\begin{lemma}\label{produkttrafo}
Let $\mu,\nu\in M^{b,+}(X)$. Then $\mu*\nu\in M_{\alpha_0}^{b,+}(X)$ if
 and only if 
$\mu,\nu\in M_{\alpha_0}^{b,+}(X)$. Moreover, if one of these conditions
 holds then
 \[R_{\alpha_0}(\mu*\nu)= R_{\alpha_0}(\mu)\bullet R_{\alpha_0}(\nu).\]
\end{lemma}

\begin{proof} If  $\mu,\nu$ have compact support, then the lemma 
is clear by Eq.~(\ref{deformconvo}).

In the general case, choose  compacta $(K_n)_{n\ge1}$  in $X$
 with $X=\bigcup_{n}K_n$ and $K_{n+1}\supset K_n$ for  $n\in\nat$.
Put 
$\mu_n:=\mu|_{K_n}$ and $\nu_n:=\nu|_{K_n}$. As the $\mu_n*\nu_n$ have
 compact support, we have
$$
\int_X \alpha_0\> d(\mu_n*\nu_n) = 
\int_X \int_X \alpha_0(x*y)\> d\mu_n(x)\> d\nu_n(y)
=\int_X \alpha_0\> d\mu_n \cdot \int_X \alpha_0\> d\mu_n.
$$
Monotone convergence implies that 
$$\int\alpha_0 \> d\mu*\nu=\int_X \alpha_0\> d\mu \cdot
 \int_X \alpha_0\> d\mu$$
where one term is finite if and only if so is  the other one. 
This proves the first part 
of the lemma.
Moreover, if these terms are finite, then the same  monotone 
convergence  argument shows that
for all $f\in C_b(X)$ with $f\ge0$,
$$\int f\> d(R_{\alpha_0}(\mu*\nu)) = 
\int f\>d( R_{\alpha_0}(\mu)\bullet R_{\alpha_0}(\nu)).$$
This implies
$R_{\alpha_0}(\mu*\nu)= R_{\alpha_0}(\mu)\bullet R_{\alpha_0}(\nu)$.
\end{proof}

\begin{remark}
Notice that the mapping $R_{\alpha_0}: M^1_{\alpha_0}(X) \to M^1(X)$ 
is not (weakly or vaguely)
continuous whenever $\alpha_0$ is unbounded. In fact, 
choose $(x_n)_{n\ge1}\subset X$ with
$\alpha_0(x_n)\to\infty$ and $\alpha_0(x_n)\ge1$.
 Then the measures
$\mu_n:= (1-\alpha_0(x_n)^{-1})\delta_e + \alpha_0(x_n)^{-1}\delta_{x_n}$
 tend to $\delta_e$ while
$$
R_{\alpha_0}(\mu_n)= 
\frac{1}{2-\alpha_0(x_n)^{-1}}((1-\alpha_0(x_n)^{-1})\delta_e +\delta_{x_n})$$
does not tend to $\delta_e= R_{\alpha_0}(\delta_e)$. 
\end{remark}

We now  investigate  convolution semigroups.

\begin{definition} A family $(\mu_t)_{t\ge0}\subset M^1(X)$ is
 called a convolution semigroup 
on $(X,*)$, if $\mu_0=\delta_e$, if $\mu_{s+t}=\mu_s*\mu_t$ for
 $s,t\ge0$, and if
the mapping $[0,\infty[\to M^1(X)$, $t\mapsto\mu_t$ is weakly 
continuous. It is well-known
(see  Rentzsch \cite{Re}) that each convolution semigroup  $(\mu_t)_{t\ge0}$ 
admits a  L\'evy measure $\eta\in M^+(X\setminus\{e\})$ 
which is characterized by
$$\int f\> d\eta= \lim_{t\to0} \frac{1}{t} \int f\> d\mu_t 
\quad \quad{\rm for}\quad \quad f\in C_c(X) 
\quad{\rm with}\quad
 e\not\in supp\> f.
$$
$(\mu_t)_{t\ge0}$ is called \emph{Gaussian}, if $\eta=0$  which 
is equivalent to saying that for all neighborhoods $U$ 
of $e\in X$, 
 $\lim_{t\to0} \frac{1}{t} \mu_t(X\setminus U)=0$.
\end{definition}

We next study  under which  conditions convolution
 semigroups on  $(X,*)$
can be deformed w.r.t. $\alpha_0$.
We here  need the following  condition on
   $\alpha_0$.

\begin{definition}
A positive  semicharacter $\alpha_0$ on $(X,*)$ is called
 \emph{exponential} if 
there exists a neighborhood $U$ of $e\in X$ and a constant
 $C>0$ such that for all $x,y\in X$ with
$y\in x*U$, $\alpha_0(y)/\alpha_0(x)\le C$.
\end{definition}

We conjecture that  positive  semicharacters are always exponential. 
Unfortunately we are not able to prove this. However,  we  present at least
some criteria and examples in Section 4 below. The following theorem  
 is motivated by \cite{Hu}, \cite{S},  where
a variant for  the group case is studied.

\begin{theorem}\label{characterization-growth}
 Let $\alpha_0$ be an  exponential positive 
 semicharacter on $(X,*)$ with $\alpha_0\ge1$.
 Then the following statements are equivalent for  a convolution semigroup
 $(\mu_t)_{t\ge0}$   on $(X,*)$ with L\'evy measure $\eta$.
\begin{enumerate}
\item[\rm{(1)}] $\mu_t\in M^1_{\alpha_0}(X)$ holds  for some $t>0$.
\item[\rm{(2)}] $\mu_t\in M^1_{\alpha_0}(X)$ holds for all  $t\ge0$, the mapping 
$\phi:[0,\infty[\to]0,\infty[$ given by $\phi(t)=\int \alpha_0\> d\mu_t$
 is continuous and multiplicative, and $(R_{\alpha_0}(\mu_t))_{t\ge0}$
 is a convolution semigroup on $(X,\bullet)$. 
\item[\rm{(3)}] For any neighborhood $U$ of
 $e\in X$, $\int_{X\setminus U} \alpha_0\> d\eta<\infty$.
\end{enumerate}
If one and hence all of these statements hold, then $\alpha_0\eta$ 
is the L\'evy measure of the convolution semigroup 
 $(R_{\alpha_0}(\mu_t))_{t\ge0}$ on $(X,\bullet)$.
\end{theorem}

In particular,  Gaussian   semigroups on  $(X,*)$ 
always lead  to Gaussian semigroups  on $(X,\bullet)$.

\begin{proof}
{\bf $(1)\Longrightarrow(2)$:} Lemma \ref{produkttrafo} implies that  $\phi\ge
1$
is well-defined and multiplicative. To check continuity,
%To check continuity, we first claim 
% that for any $r<1$ there exists $t_0>0$ such that for
% $0\le s\le t_0$, $\phi(s)\ge r$ holds. This follows from the
% fact that $U:=\{x\in X:\> \phi(x)\ge\sqrt r\}$ is a neighborhood of $e$,
% and that $\mu_s(U)\ge \sqrt r$ for $s$ sufficiently small.
we observe that the multiplicativity  implies that 
 for $N\in\nat$  and $0\le s\le 1/N$, 
 $\phi(s)^{N}\phi(1-sN)=\phi(1)$  and hence $\phi(s)\le \phi(1)^{1/N}\to 1$
 for $N\to\infty$.
Therefore,  $\phi$ is continuous at $t=0$ and hence,
 as a multiplicative function, on $[0,\infty[$.
Using Lemma  \ref{produkttrafo} and the fact that the mapping
$[0,\infty[\to M^1(X)$, $t\mapsto  R_{\alpha_0}(\mu_t)= \phi(t)^{-1} \alpha_0\mu$ is vaguely
and  hence  weakly continuous, we conclude that 
 $(R_{\alpha_0}(\mu_t))_{t\ge0}$ is a
  convolution semigroup on $(X,\bullet)$.

{\bf $(2)\Longrightarrow(3)$:} The measure $\rho:={\bf 1}_{\{\alpha_0\ge 2\}}\eta\in M^{b,+}(X)$
is the L\'evy measure of the Poisson semigroup
 $(\nu_t:=e^{-\|\rho\|t}\cdot exp(t\rho))_{t\ge0}$, 
$exp$ denoting the exponential function on the Banach 
algebra $(M_b(X),*)$. Moreover, it is easy to see that
  $\eta-\rho$ is the  L\'evy measure of 
 a further convolution semigroup $(\tilde \nu_t)_{t\ge0}$ with
 $\mu_t=\nu_t*\tilde\nu_t$ for $t\ge0$. Lemma  \ref{produkttrafo} shows that
  $\nu_t\in M^1_{\alpha_0}(X)$ for $t\ge0$. As obviously
 $\rho\le (e^{\|\rho\|t}/t)\nu_t$ for $t>0$, we obtain
$\rho\in M_{\alpha_0}^{b,+}(X)$ and  thus (3).
Furthermore, for $f\in C_c(X)$ with $e\not\in supp\> f$, 
$$\lim_{t\to 0}\frac{1}{t}\int f\> dR_{\alpha_0}(\mu_t) =
\lim_{t\to 0}\frac{1}{t}\int f\alpha_0\> d\mu_t= \int  f\alpha_0\> d\eta.$$
Hence, $\alpha_0\eta$ is the  L\'evy measure of the 
 semigroup  $(R_{\alpha_0}(\mu_t))_{t\ge0}$ on $(X,\bullet)$.
\end{proof}

The proof of  $(3)\Longrightarrow(1)$ is more involved.  Recapitulate that  
for a convolution semigroup $(\mu_t)_{t\ge0}$ on $(X,*)$,
the translation operators $T_t(f):=\mu_t^-*f$ ($t\ge0$) form
 a strongly continuous, positive contraction
 semigroup  on $L^1(X,\omega)$,
 $\omega$ being the Haar measure of $(X,*)$; see [BH]. Let $A$ be its infinitesimal generator 
 with the dense domain $D_A\subset L^1(X,\omega)$. We  have:

\begin{lemma}\label{nettererzeuger} Let $\alpha_0$ be a  positive  semicharacter 
 and  $(\mu_t)_{t\ge0}$   a convolution semigroup  on $(X,*)$  whose 
L\'evy measure $\eta$ satisfies   $\int_{\{\alpha_0\ge2\}}\alpha_0\> d\eta<\infty$.
 Then for each 
neighborhood $U$ of $e\in X$  there exists
$f\in C_c(X)\cap D_A$ with $f\ge0$, $f=f^*\ne0$, $supp\> f\subset U$,
 and $\int |Af|\alpha_0\> d\omega<\infty$.
\end{lemma}

 \begin{proof}
Let $U$ be a compact neighborhood of $e\in X$ with $U^-=U$.
 Then by \cite{Re}, there exists $f\in D_A$ with $\int_X 
f\>d\omega=1$, $f\ge0$, $f=f^*$,
 and $supp\> f\subset U$.
 Let $x\not\in U*U$ and $y\in U$. Then $f(x*y)=0$, which means that  the translate $f_x$ given by $f_x(y):=f(x*y)$ satisfies $f_x=0$ on $U$, and hence
$$Af(x)=\lim_{t\to0} \frac{1}{t} (\mu_t^-*f_x(e)- f_x(e))=\int f(x*y)\> d\eta(y).$$
Consequently, by Fubini's theorem,
\begin{align}
\int_X |Af|\cdot \alpha_0\> d\omega 
&= \int_{U*U}|Af|\cdot \alpha_0\> d\omega \>+\>
\int_{X\setminus U*U}|Af|\cdot \alpha_0\> d\omega 
\notag\\
&\le \int_{U*U}|Af|\cdot \alpha_0\> d\omega \>+\>
 \int_X\int_{X\setminus U*U} f(x*y)\alpha_0(x)\> d\omega(x)\> d\eta(y).
\notag
\end{align}
Now
\begin{align}
\int_{X\setminus U*U} f(x*y)\alpha_0(x)\> d\omega(x) &=
\int_X f(x) ({\bf 1}_{X\setminus U*U}\alpha_0)(x*\bar y) d\omega(x)
\notag\\
&\le {\bf 1}_{X\setminus U}(\bar y)\cdot
 \int_X f(x)\alpha_0(x)\alpha_0(\bar y) d\omega(x)
\notag\\
&
=  {\bf 1}_{X\setminus U}(y)\alpha_0(y)\cdot  \int f\alpha_0 d\omega.
\notag
\end{align}
As $\int {\bf 1}_{X\setminus U}\alpha_0\> d\eta<\infty$ by assumption,
 $\int |Af|\alpha_0\>d\omega <\infty$ as claimed.
\end{proof}

 $(3)\Longrightarrow(1)$ in the theorem now follows from Lemma
\ref{nettererzeuger} and the following result.

\begin{lemma} Let $\alpha_0$ be an exponential  positive  semicharacter with
 $\alpha_0\ge1$,
 and  $(\mu_t)_{t\ge0}$   a convolution semigroup  on $(X,*)$ with generator $A$. Assume that for each 
neighborhood $U$ of $e\in X$   there exists
$f\in C_c(X)\cap D_A$ with $f\ge0$, $f=f^*\ne0$, $supp\> f\subset U$
 and $\int |Af|\alpha_0\> d\omega<\infty$.
Then for all $t\ge0$, $\int \alpha_0\> d\mu_t<\infty$.
\end{lemma}

\begin{proof}
Let $U$ be a neighborhood of $e\in X$ and $C_1>0$ a constant with
 $C_1\alpha_0(x)\le\alpha_0(z)$ for  $x\in X$ and $z\in U*x$. Let
$f\in C_c(X)\cap D_A$ with $f\ge0$, $f=f^*\ne0$, $supp\> f\subset U$
 and $\int |Af|\alpha_0\> d\omega<\infty$. Then for all $m\in\nat$,
 the functions $\alpha_m:=\alpha_0\wedge m\in C_b(X)$ also satisfy 
 $C_1\alpha_m(x)\le\alpha_m(z)$ for  $x\in X$, $z\in U*x$.
Hence, there is a constant $C_2>0$ depending on $f$ such that for all  $m\in\nat$ and
 $x\in X$,
\begin{equation}\label{absch1}
\alpha_m(x)\le C_2\cdot\int\alpha_m(x*y)f(y)\> d\omega(y)=C_2\cdot \alpha_m*f(x).
\end{equation}
Moreover, as $\alpha_0\ge1$, we have for all $m\in\nat$ and  $x,y\in X$,
 \begin{equation}\label{absch2}
\alpha_m(x*y)\le m\wedge\alpha_0(x*y) =  m\wedge(\alpha_0(x)\alpha_0(y))\le \alpha_m(x)
 \alpha_m(y).
\end{equation}
Define $h_m(t):= \int (\mu_t*f )\cdot\alpha_m\>d\omega= \int\alpha_m*f  d\mu_t$.  As
$f\in D_A$ and $Af\in L^1(X,\omega)$ holds, we obtain
$\frac{d}{dt}\mu_t*f=\mu_t*Af$ and hence
$$h_m^\prime(t)= \int (\mu_t*Af)\cdot  \alpha_m\>d\omega= 
 \int  \int \alpha_m(x*y) Af(y) \> d\mu_t(x)\>d\omega(y) .$$
Therefore, by (\ref{absch2}) and  (\ref{absch1}),
\begin{align}
|h_m^\prime(t)|&\le \int  \int \alpha_m(x*y) |Af(y)| \> d\mu_t(x)\>d\omega(y)
\>\le \> \int\alpha_m  d\mu_t \cdot  \int\alpha_m |Af|\> d\omega
\notag\\
&\le C_2  \int\alpha_m*f  d\mu_t \cdot  \int\alpha_m |Af|\> d\omega
\>\le \>  C_2  \int\alpha_0 |Af|\> d\omega\cdot h_m(t).
\notag
\end{align}
This yields $h_m(t)\le h_m(0)\> e^{tC}$ 
for $t\ge 0$ and some constant $C\ge0$ independent of $m$. Hence,
again by  (\ref{absch1}),
$$\int \alpha_m\>d\mu_t\le  C_2  \int\alpha_m*f  d\mu_t=  C_2 h_m(t)\le
 c_2 e^{tC} \int\alpha_0 f\>d\omega$$
for all $m\in\nat$.  This yields the claim $\int \alpha_0\>d\mu_t<\infty$ for $t\ge0$.
\end{proof}

Notice that the growth condition on $\alpha_0$ was  needed above only for 
 the  preceding  lemma. 
Theorem \ref{characterization-growth} therefore admits the following variant.

\begin{theorem}
Let $\alpha_0$ be a positive semicharacter and $(\mu_t)_{t\ge0}\subset M^1(X)$
a Poisson semigroup on $(X,*)$, which means that
$\mu_t=e^{-t\|\rho\|}exp(t\rho)$ for all $t\ge0$ and some $\rho\in
M^{b,+}(X)$. Then $\rho$ is the L\'evy measure of  $(\mu_t)_{t\ge0}$, and the
statements (1)--(3) of Theorem \ref{characterization-growth} are equivalent.
\end{theorem}

\begin{proof} It suffices to check $(3)\Longrightarrow(1)$. However, if
  $R:=\int\alpha_0\> d\rho<\infty$, then 
for all $n\ge 0$, $\int\alpha_0\> d\rho^{(n)}=R^n$
and hence $\int\alpha_0\> d\mu_t<\infty$ for all $t\ge0$.
\end{proof}

\begin{remark}
Let  $\alpha_0$ be an exponential positive
 semicharacter and $(\mu_t)_{t\ge0}\subset  M^1_{\alpha_0}(X)$  
a convolution semigroup  on $(X,*)$.
Then the convolution operators $(T_t)_{t\ge0}$  on $C_0(X)$ with
 $T_tf:=\mu_t^-*f$ form a Feller semigroup. Its  generator $A$ 
with
$$Af(x)= \lim_{t\to 0} \frac{1}{t}(\mu_t^-*f(x) -f(x)) \quad\quad
(x\in X, \> f\in D(A))$$
admits a 
 $\|.\|_\infty$-dense  domain $D(A)$  in $C_0(X)$; see  \cite{ReV}.
Now consider the generator $A^{\alpha_0}$ of the
 Feller semigroup on $C_0(X)$ which is associated with the renormalized
 convolution semigroup $(R_{\alpha_0}(\mu_t))_{t\ge0}$ on $(X,\bullet)$. 
Using the notation above, we have
$$((R_{\alpha_0}\mu_t)^-\bullet f)(x)= \frac{1}{\phi(t)} ((\alpha_0\mu_t)^-
\bullet f)(x)=  \frac{1}{\phi(t)\alpha_0(x)}(\mu_t*\alpha_0f)(x).$$
Theorem 
\ref{characterization-growth}(2) shows that
 $\phi(t)=e^{ct}$ for
some $c\in\reel$, and 
$$\lim_{t\to0}  \frac{1}{t}(1/\phi(t)-1)=-c.$$
 Hence
\begin{align}
A^{\alpha_0}f(x)&=\lim_{t\to0}  \frac{1}{t}\Bigl(
 \frac{1}{ \phi(t)\alpha_0(x)}(\mu_t*\alpha_0f)(x) -f(x)\Bigr) \notag\\
&=  \frac{1}{ \alpha_0(x)}\lim_{t\to0}  \frac{1}{t}\Bigl(
 \frac{1}{ \phi(t)} (\mu_t*\alpha_0f)(x) -(\alpha_0f)(x)\Bigr)
 \notag\\
&= \frac{1}{ \alpha_0(x)} A(\alpha_0f)(x) + \frac{1}{ \alpha_0(x)}
\lim_{t\to0}\Bigl(  \frac{1}{t}(1/h(t)-1) (\mu_t*\alpha_0f)(x)\Bigr)
 \notag\\
&= \frac{1}{ \alpha_0(x)} A(\alpha_0f)(x) -c f(x). \notag
\end{align}
Therefore 
\begin{equation}
A^{\alpha_0}= M_{1/\alpha_0}\circ A\circ  M_{\alpha_0} -c
\end{equation}
at least on $D(A^{\alpha_0})\cap C_c(X)$,
where $M_g$ denotes the multiplication operator with 
 $g\in C(K)$. The same holds for other function spaces like $L^p(X,\omega)$. 
\end{remark}

\section{Exponential positive semicharacters}

It seems reasonable to  conjecture that  positive  semicharacters are always 
exponential. Unfortunately we are not able to prove this.
Here are, at least, 
some criteria and  several examples:

\begin{lemma}
\begin{enumerate}
\item[\rm{(1)}] If $(X,*)$ is discrete, then 
 $\alpha_0$ is always
  exponential.
\item[\rm{(2)}] Let $\alpha_0,\alpha_1$ be exponential
 positive  semicharacters on  $(X,*)$, and let $(X,\bullet)$ be the
deformation of  $(X,*)$ w.r.t.  $\alpha_0$. Then $\alpha_1/\alpha_0$ is
an exponential positive  semicharacter on $(X,\bullet)$.
\end{enumerate}
\end{lemma}

\begin{proof} 
Part (1) is clear by taking $U=\{e\}$. For the proof of (2) choose
neighborhoods $U_0,U_1$ of $e$ and constants $C_0,C_1$ associated with 
$\alpha_0,\alpha_1$ respectively. For $U:=U_0\cap U_1\cap U_0^-\cap U_1^-$ and
$C:=C_0C_1$, we obtain that for $x,y\in X$ with
$y\in x*U$, we have $x\in y*U$ and thus 
$\alpha_0(x)\alpha_1(y)/(\alpha_0(y)\alpha_1(x))\le C$ as claimed.
\end{proof} 

% Here are a few   examples:

\begin{example}
 In \cite{Z2}, Zeuner presented  quite general, but technical
  conditions on a function $A\in C([0,\infty[)\cap C^1(]0,\infty[)$ with
  $A(x)>0$ for $x\ge0$ which ensures that there exists a unique commutative
  hypergroup $([0,\infty[,*)$ whose semicharacters are precisely
  the eigenfunctions of the Sturm-Liouville operator
$$L_Af:= -f^{\prime\prime}-(A^\prime/A)f^\prime $$
with  initial conditions $f(0)=1$ and $f^\prime(0)=0$; see also Section 3.5
of \cite{BH}. This hypergroup 
is called the   Sturm-Liouville hypergroup associated with $A$. Moreover, to
the knowledge of the authors, all known hypergroup structures on $[0,\infty[$
appear in this way (up to isomorphism); see also \cite{BH} for details.
We here mention that Zeuner's approach in particular includes all 
Chebli-Trimeche hypergroups and thus all double coset
 hypergroups associated with noncompact symmetric spaces of rank one.

We claim that all positive semicharacters on a   Sturm-Liouville hypergroup on
$[0,\infty[$ with $A$ satisfying Zeuner's conditions are exponential.
 To prove this, recall from Section 3.5 in  \cite{BH}
that Zeuner's conditions imply that
\begin{equation}\label{rholimit}
\rho:=\frac{1}{2}\lim_{x\to\infty} A^\prime(x)/A(x) \ge0
\end{equation}
exists, and that the  positive semicharacters are precisely  the
unique solutions $\phi_{i\lambda}$ of
$$L_A\phi_{i\lambda}=(\rho^2-\lambda^2)\phi_{i\lambda},\quad\quad
\phi_{i\lambda}(0)=1,\> \phi_{i\lambda}^\prime(0)=0$$
with $\lambda\ge 0$. Moreover, the renormalization 
$([0,\infty[,\bullet)$ of  $([0,\infty[,*)$ w.r.t.  $\phi_{i\lambda}$ is again
a  Sturm-Liouville hypergroup associated with the renormalized function 
 $A_\lambda:=\phi_{i\lambda}^2\cdot A$ where  $A_\lambda$ again satisfies 
 Zeuner's conditions; see Section 3.5.51 of  \cite{BH}. Applying
 (\ref{rholimit}) to
 $A$ as well as to  $A_\lambda$, we see that $\lim_{x\to\infty}
 \phi_{i\lambda}^\prime(x)/ \phi_{i\lambda}(x)$ exists. As 
  $supp(\delta_x*\delta_y)\subset [|x-y|, x+y]$  
 for $x,y\ge 0$, it follows from the mean-value theorem 
 that $\phi_{i\lambda}$ is
exponential.
\end{example}

\begin{example}
Let $V$ be a finite-dimensional Euclidean vector space, $K\subset O(V)$ a
compact subgroup, and $V^K$ the associated orbit hypergroup as in Example
\ref{ex_1}(2). Then, for each $\rho\in V$, the positive semicharacter
$\alpha_{i\rho}$ with
$\alpha_{i\rho}(K.x)=\int_K e^{-\langle \rho, k.x\rangle}dk$ ($x\in V$) is 
exponential.
In fact, we may take $U:=\{K.x:\> x\in V, \> \|x\|_2\le1\}\subset V^K$ as a
neighborhood of the identity. For orbits $K.x,K.y\in V^K$ with $K.x\in U*K.y$
we then have representatives $x,y\in V$ with $\|x-y\|_2\le 1$ which 
implies that 
$e^{-\langle \rho, k.x\rangle}\le
e^{-\langle \rho, k.y\rangle} e^{\|\rho\|_2}$  for $k\in K$ and thus
$\alpha_{i\rho}(K.x)\le \alpha_{i\rho}(K.y)e^{\|\rho\|_2}$ as claimed.
\end{example}

\begin{example}
Let $G$ be a (not necessarily complex) noncompact, connected
semisimple  Lie group with finite center
and $K$  a maximal compact subgroup.
Let $G=NAK$ and  $\mathfrak g= \mathfrak n +\mathfrak a+ \mathfrak k$
be the corresponding Iwasawa decompositions. For $g\in G$ let
 $A(g)\in \mathfrak a$  be the unique element with $g\in N \> exp(A(g))K$.
Let  $\Sigma^+$ be the set of  positive roots (for the order corresponding to 
$\mathfrak n$), and 
$\rho=\frac{1}{2}\sum_{\alpha\in\Sigma^+} m_\alpha\alpha$ 
the half sum of positive  roots with $m_\alpha$ as multiplicity of 
$\alpha$. 
 Then, by a formula of Harish-Chandra
 (see Theorem IV.4.3 of \cite{He}), the spherical functions on $G$, i.e., 
 the multiplicative functions on $G//K$,
 are given by
$$\phi_\lambda(g)=\int_K e^{\langle i\lambda+\rho, A(kg)\rangle} \> dk 
\quad\quad (g\in G),$$
where $\lambda$ runs through $\mathfrak a_\comp$, the complexification of 
$\mathfrak a$.
Clearly, the $\phi_\lambda$ for $\lambda\in i\cdot \mathfrak a$  are positive 
multiplicative functions. These functions are also exponential.
To prove this, we conclude from   Lemma  IV.4.4 of \cite{He} that
$$\phi_\lambda(g^{-1}h)=\int_K e^{\langle- i\lambda+\rho, A(kg)\rangle}
 e^{\langle i\lambda+\rho, A(kh)\rangle} \> dk 
\quad\quad (g,h\in G).$$
Hence, for each compact neighborhood $U$ of $e$ there is a constant $C>0$
 such that $\phi_\lambda(g^{-1}h)\le C\phi_\lambda(h)$ for all $g\in U$, $h\in
 G$ and  $\lambda\in i\cdot \mathfrak a$.
\end{example}

\begin{example}
Let $R$ be a (reduced, not necessaryly crystallographic)
 root system in $\reel^n$ with the standard  inner
 product $\langle\,.\,,\,.\,\rangle,$  i.e.
 $R\subset \reel^n\setminus\{0\}$ is finite with $R\cap \reel\alpha =
 \{\pm\alpha\}$  and $\sigma_\alpha(R) = R$ for all $\alpha\in R$, where $\sigma_\alpha$ is
 the reflection in the hyperplane perpendicular to $\alpha$. Assume also
without loss of generality for our considerations that $\langle\alpha,\alpha\rangle = 2$ for all
$\alpha\in R$. Let $W$ be the finite reflection group generated by 
the $\sigma_\alpha$
and let
$k: R\to [0,\infty[$  be a fixed
multiplicity function on $R,$  i.e. a function
which is constant on the orbits under the action of $W$.
The   (so-called rational) Dunkl operators attached to
$G$ and $k$ are defined by
\begin{equation}\label{(1.10a)}
 T_\xi(k) f(x)  \,=\, \partial_\xi f(x)  +
\sum_{\alpha\in R_+} k(\alpha) \langle \alpha,\xi\rangle
\frac{f(x) - f(\sigma_\alpha x)}{\langle\alpha,x\rangle},\quad x,\,\xi\in
 \reel^n.
\end{equation}
Here $\partial_\xi$ denotes the derivative in direction $\xi$ and $R_+$ is
 some fixed positive subsystem  of $R$. The definition is independent of the
 special choice of $R_+$, due to the $G$-invariance of $k$.
 As first shown in \cite{D1}, the
$T_\xi(k), \,\xi\in \reel^n$
 generate a commutative algebra of differential-reflection operators. 
 This is the foundation for  rich analytic structures related
with them.
In particular, there  exists a 
  counterpart of the  exponential function, 
the Dunkl kernel, and an analogue of the Euclidean Fourier transform
with respect to this kernel. The Dunkl kernel $E_k$ is  holomorphic
on $\comp^n\times \comp^n$ and symmetric in its arguments. 
  Similar to
spherical functions
on  symmetric spaces, the function $E_k(\,.\,,y)$ with fixed
 $y\in  \comp^n$ may  be
characterized as  unique analytic solution of the
 joint eigenvalue problem
\begin{equation}\label{(1.104)}
 T_\xi(k) f\,=\, \langle\xi,y\rangle f 
\quad\text{for all }\,\xi \in  \comp^n,\,\,\,
f(0) = 1;
\end{equation}
 c.f. \cite{O}.
Apart from the trivial case $k=0$ with $E_k(x,y) = e^{\langle x,y\rangle}$,
 $E_k$  is explicitly known in a few cases only like  $n=1$; see  
\cite{Ro3} for a survey.
The $G$-invariant counterpart of $E_k$ is the 
generalized Bessel function
\[ J_k(x,y) = \frac{1}{|G|} \sum_{g\in G} E_k(gx,y)\]
 which  is $G$-invariant in $x, y$ 
and naturally considered on the closed  positive Weyl chamber $C$  associated
with $R_+$. For $n=1$,  $J_k$ is 
 a usual Bessel function. Moreover, 
in the cristallographic case and for certains
half-integer multiplicities,
the $J_k$ are  the multiplicative functions of certain  
Euclidean orbit hypergroups as in Example \ref{ex_1}.
Here, and for  $n=1$,   the  $J_k(x,y)$ ($y\in\comp^n$) therefore
form the multiplicative functions of
 some commutative hypergroup  on $C$.
It is conjectured that there exist such
 commutative hypergroups on $C$ for all root systems and 
multiplicities $k\ge0$. 
Only part of this conjecture has been verified
up to now in \cite{Ro2}.

Now  fix a root system $R$ and $k\ge 0$ such that 
the  $J_k(.,y)$ ($y\in\comp^n$) are  the multiplicative functions
of a  commutative hypergroup  $(C,*)$. To find positive semicharacters, we
employ 
the following psoitive integral representation for  $E_k$ (and thus
 $J_k$):
For given $R$, $k\ge 0$, and $x\in\reel^n$ there exists a unique pribability
measure $\mu_x$ on $\reel^n$ such that
\begin{equation}\label{pos-int}
J_k(x,y)=\int e^{\langle z,y\rangle}\> d\mu_x(z)
\quad\quad{\rm for}\quad y\in\comp^n.
\end{equation}
Moreover,  
$supp\>\mu_x\subset \{ z\in\reel^n: \> \|z\|_2\le \|x\|_2\}$.
Thus, for each $y\in\reel^n$,  $J_k(.,y)$ is a positive semicharacter on
$(C,*)$. We claim that these semicharacters are exponential.

To show this, let $U:=\{z\in C:\> \|z\|_2\le1\}$ and $x_1,x_2\in C$ with
$x_1\in U*x_2$. We conclude from Theorem 4.1 of \cite{Ro2} that then 
$x_1\in C\cap\bigcap_{w\in W}\{z\in\reel^n:\> |z-w.x_2|\le1\}$ holds.
As $\|z-w\|\le\|z-w.x\|$ for all $x,z\in C$ and $w\in W$ by Ch. 3 of 
\cite{GB}, we even have $\|x_1-x_2\|\le 1$. In the same way as in Example 4.3
we now obtain from Eq.~(\ref{pos-int}) that $J_k(x_1,y)\le
e^{\|y\|}J_k(x_2,y)$
which proves that $J_k(.,y)$ is exponential for each $y\in \reel^n$. 
\end{example}


\begin{thebibliography}{AA}

\bibitem{BH} W. Bloom, H. Heyer, \textit{Harmonic Analysis of 
    Probability Measures on Hypergroups.} De Gruyter-Verlag, Berlin, 1994.

\bibitem{D1} C.F. Dunkl, \textit{Differential-difference operators
associated to reflection groups}, Trans. Amer. Math. Soc. \textbf{311} (1989),
 167--183.

\bibitem{GB} L.C. Grove, C.T. Benson, \textit{Finite Reflection Groups}.
Springer-Verlag, 2nd ed., 1985.

\bibitem{He} S. Helgason, \textit{Groups and Geometric Analysis.} American
Mathematical Society, 2000.

\bibitem{Hu} A. Hulanicki,  \textit{A class of convolution semi-groups of 
measures on a Lie group}. In: Probab. Theory on Vector Spaces II, Proceedings.
 Lecture Notes in Mathematics Vol. 828, Springer 1980, pp. 82--101. 


\bibitem{J} R.I. Jewett, \textit{Spaces with an abstract convolution
of measures}. Adv. Math. \textbf{18} (1975), 1--101.

\bibitem{Kl1} A.  Klyachko, \textit{Stable bundles, representation theory, and
    Hermitian operators}, Selecta Math. (new Series)  \textbf{4} (1998),
419--445. 

\bibitem{Kl2} A.  Klyachko, \textit{Random walks on symmetric spaces
and inequalities for matrix spectra}, Linear Algebra Appl. 
\textbf{319} (2000), 37--59.

\bibitem{Kn} A.W.  Knapp, \textit{Lie Groups beyond an Introduction.}
 Birkh\"auser, Boston, 1996.

\bibitem{KT} A. Knutson, T. Tao,  \textit{The honeycomb model of $GL_n(\comp)$
tensor product I. Proof of the saturation conjecture}, J. Amer. Math. Soc. 
\textbf{12} (1999), 1055-1090.


%\bibitem{Ko} T. Koornwinder, \textit{Jacobi functions and 
%analysis of noncompact semisimple
%Lie groups}. In: R. Askey et al. (eds.), Special Functions:
%Group Theoretical Aspects
%and Applications. Dordrecht-Boston-Lancaster: D. Reidel
%Publishing Company 1984, pp. 1-85.

\bibitem{O} E.M. Opdam, \textit{Dunkl operators, Bessel functions and the
discriminant of a finite Coxeter group}, Compos. Math.
\textbf{85} (1993), 333--373.

\bibitem{Re} C. Rentzsch, \textit{A L\'evy-Khintchine type representation 
of convolution semigroups on commutative hypergroups.}
Probab. Math. Stat. \textbf{18}  (1998), 185-198. 

\bibitem {ReV} C.~Rentzsch, M.~Voit, \textit{ L\'evy processes on
 commutative hypergroups}. Contemp. Math.  \textbf{261} (2000), 83--105. 

\bibitem {RW} L.C.G.~Rogers, D.~Williams,
 \textit{  Diffusions, Markov Processes,
 and Martingales, Vol. II}. Wiley:
Chichester -- New York, 1987.

\bibitem{Ro1} M. R\"osler, \textit{Positivity of Dunkl's
 intertwining operator.}
Duke Math. J. \textbf{98}  (1999),  445-463.


\bibitem{Ro2} M. R\"osler, \textit{A positive radial product 
formula for the Dunkl kernel}. Trans. Amer. Math. Soc.
\textbf{355} (2003), 2413-2438.

\bibitem{Ro3} M. R\"osler, \textit{Dunkl operators: Theory and Applications.}
 Lecture Notes for the SIAM Euro Summer School: Orthogonal Polynomials 
and Special Functions. (Leuven,  2002). 
In: Springer Lecture Notes in Mathematics, Vol. 1817, 2003, pp. 93 - 136. 

\bibitem{RV} M. R\"osler, M. Voit, \textit{$SU(d)$-biinvariant random walks on
    $SL(d,\comp)$ and their Euclidean counterparts}, preprint 2003.

\bibitem{S} E.~Siebert,  \textit{Continuous convolution semigroups 
integrating a submultiplicative function}, Man. Math. \textbf{37} (1982),
 383--391.

\bibitem{V1} M. Voit, \textit{ Positive
 characters on commutative
hypergroups and some  applications}, Math. Z. \textbf{198} (1988), 405-421.

\bibitem{V2} M. Voit,
 \textit{A Girsanov-type formula for L\'evy processes 
on commutative hypergroups.} 
In: H. Heyer et al. (eds.), Infinite Dimensional Harmonic Analysis,
 Proc. Conf.,
Gr\"abner 2000, pp. 346--359.


\bibitem{Z1} H. Zeuner, \textit{
The central limit theorem for Ch\'ebli-Trimeche hypergroups.} 
J. Theor. Probab. \textbf{2} (1989), 51-63.

\bibitem{Z2} H. Zeuner, \textit{ 
Moment functions and laws of large numbers on hypergroups.} 
Math. Z. \textbf{211} (1992),  369-407. 

\end{thebibliography}
\end{document}